\documentclass{article}
\usepackage{mathtools}

\usepackage{arxiv}

\usepackage[utf8]{inputenc} % allow utf-8 input
\usepackage[T1]{fontenc}    % use 8-bit T1 fonts
\usepackage{hyperref}       % hyperlinks
\usepackage{url}            % simple URL typesetting
\usepackage{booktabs}       % professional-quality tables
\usepackage{amsfonts}       % blackboard math symbols
\usepackage{nicefrac}       % compact symbols for 1/2, etc.
\usepackage{microtype}      % microtypography
\usepackage{lipsum}
\usepackage{graphicx}
\usepackage{floatrow}
\usepackage{graphicx}
\usepackage[label font=bf, 
            labelformat=simple]{subfig}
\usepackage{caption}
\usepackage{lineno}
\usepackage{svg}
\usepackage{parskip} % For paragraph layout
\usepackage{setspace} % For using single or double spacing
\usepackage{emptypage} % To insert empty pages
\usepackage{multicol} % To write in multiple columns (executive summary)
\usepackage{epstopdf}
\PrependGraphicsExtensions{.tif, .tiff}
\usepackage{titlesec}
\titlespacing{\section}{0pt}{3.3ex}{2ex}
\titlespacing{\subsection}{0pt}{3.3ex}{1.65ex}
\titlespacing{\subsubsection}{0pt}{3.3ex}{1ex}
\usepackage{color}
\usepackage{gensymb}
\usepackage[english]{babel} % The document is in English  
\usepackage[utf8]{inputenc} % UTF8 encoding
\usepackage[T1]{fontenc} % Font encoding
\usepackage[11pt]{moresize} % Big fonts

\usepackage{graphicx}
\usepackage{transparent} % Enables transparent images
\usepackage{eso-pic} % For the background picture on the title page
\usepackage{subfig} % Numbered and caption subfigures using \subfloat.
\usepackage{tikz} % A package for high-quality hand-made figures.
\usetikzlibrary{}
\graphicspath{{./Images/}} % Directory of the images
\usepackage{caption} % Coloured captions
\usepackage{xcolor} % Coloured captions
\usepackage{amsthm,thmtools,xcolor} % Coloured "Theorem"
\usepackage{float}

\usepackage{amsmath}
\usepackage{amsthm}
\usepackage{amssymb}
\usepackage{amsfonts}
\usepackage{bm}
\usepackage[overload]{empheq} % For braced-style systems of equations.
\usepackage{fix-cm} % To override original LaTeX restrictions on sizes

\usepackage{tabularx}
\usepackage{longtable} % Tables that can span several pages
\usepackage{colortbl}

\usepackage{algorithm}
\usepackage{algorithmic}

\usepackage{pdfpages} % To include a pdf file
\usepackage{afterpage}
\usepackage{lipsum} % DUMMY PACKAGE
\usepackage{fancyhdr} % For the zheaders
\usepackage{amsthm}   % for theorem, lemma, proof environments
\usepackage{amsmath}  % you probably already have this for equations

\usepackage{comment}
\graphicspath{ {./images/} }
\title{Why Is Nicomachus' Identity Special?}

\author{
Guglielmo Vesco\\
Mathematics and Physics Teacher\\
IIS P. Frisi, Milan, Italy\\
\texttt{guglielmo.vesco@gmail.com}
}

\begin{document}
\maketitle

\begin{abstract}
Nicomachus' identity

\[
\sum_{k=1}^{n}k^3
=
\left(\sum_{k=1}^{n}k\right)^2
\]

is familiar, but the structure behind it is less often emphasized. We begin with Nicomachus' original pattern, in which successive cubes are expressed as successive blocks of consecutive odd integers, and show how triangular numbers naturally enter the resulting sum-of-cubes identity. We then reverse the usual triangular-number argument and start from a general difference of squares representing $n^p$ as a sum of $n$ consecutive odd integers. This construction is possible for every $p\geq2$, but only for $p=3$ do the two squares correspond to consecutive triangular numbers, which explains why Nicomachus' odd-number blocks fit together without gaps or repetitions. Finally, using the leading terms of the Faulhaber polynomials, we show that, apart from the trivial case, Nicomachus' identity is the unique relation of the form $S_s(N)=S_1(N)^r$ for positive integer exponents.
\end{abstract}


\section{Introduction}

The identity

\begin{equation}
\sum_{k=1}^{n}k^3
=
\left(\sum_{k=1}^{n}k\right)^2
\end{equation}

is often encountered as a striking piece of recreational mathematics: the sum of the first $n$ cubes is exactly the square of the sum of the first $n$ positive integers. Its compact form naturally raises two questions. Why do cubes appear on the left, and why does a square appear on the right? More broadly, is this identity an isolated coincidence, or does it reflect a more rigid underlying structure?\\
The exceptional character of this identity has been noted before. In particular,
Edmonds investigated the special relation between sums of powers that underlies
Nicomachus' identity \cite{Edmonds1957}.

The identity is traditionally associated with Nicomachus of Gerasa, whose original observation was not stated in this modern form. Instead, Nicomachus described a pattern in which successive cubes are represented by successive blocks of consecutive odd integers. We first return to that observation and show how the familiar identity emerges from it through the elementary fact that the sum of the first $m$ odd integers is $m^2$. This naturally introduces the triangular numbers and leads to the standard telescoping relation

\begin{equation}
T_n^2-T_{n-1}^2=n^3.
\end{equation}

Modern proofs often begin with this relation. Here we reverse the point of view and ask why consecutive triangular numbers should appear in the first place. Starting from a general difference of squares,

\begin{equation}
n^p=b^2-a^2,
\end{equation}

we require it to represent a block of exactly $n$ consecutive odd integers, as in Nicomachus' pattern. The resulting construction works for every power $p\geq2$, but only in the cubic case do the two endpoints become consecutive triangular numbers. This explains, within the same framework, why the odd-number blocks for the cubes fit together without gaps or repetitions.

Finally, we ask whether the same exceptional behavior appears at the level of the global power-sum identity. Using the leading terms of the Faulhaber polynomials, we determine all positive integer pairs $(s,r)$ for which

\begin{equation}
\sum_{k=1}^{N}k^s
=
\left(\sum_{k=1}^{N}k\right)^r
\end{equation}

can hold for every positive integer $N$. Apart from the trivial case, the only possibility is precisely Nicomachus' pair $(s,r)=(3,2)$. Thus the cubic identity is singled out both locally, through its consecutive-odd-number structure, and globally, among identities of this power-sum form.

\section{Nicomachus, Odd Numbers, and Triangular Numbers}

Little is known about the life of Nicomachus of Gerasa. He probably lived around AD~100 and belonged to the Pythagorean tradition, although the dates of his life can only be estimated from the authors whom he mentions and from later writers who refer to him \cite{Dillon1996}. Near the end of Book~II, Chapter~20 of his \emph{Introduction to Arithmetic}, Nicomachus recorded the following pattern among the positive odd integers \cite[Book~II, Chapter~20]{Nicomachus1926}:

\begin{equation}
1=1^3,
\qquad
3+5=2^3,
\qquad
7+9+11=3^3,
\qquad
13+15+17+19=4^3,
\qquad\ldots
\end{equation}

Thus the first cube is represented by the first odd integer, the second cube by the next two odd integers, the third cube by the next three, and so on. In general, the $k$th cube is represented by a block of $k$ consecutive odd integers. Nicomachus stated this regularity, but did not provide a proof in the modern deductive sense.

Taking this observed pattern as our starting point, the familiar identity for the sum of cubes follows naturally. After the first $n$ blocks have been used, the total number of odd integers appearing is

\begin{equation}
1+2+\cdots+n.
\end{equation}

This is the $n$th triangular number,

\begin{equation}
T_n=1+2+\cdots+n=\frac{n(n+1)}{2}.
\end{equation}

On the other hand, the sum of the first $m$ positive odd integers is

\begin{equation}
1+3+5+\cdots+(2m-1)=m^2
\end{equation}

\cite{Sangwin2023}. Since Nicomachus' first $n$ blocks contain exactly the first $T_n$ odd integers, their total is therefore

\begin{equation}
\sum_{k=1}^{n}k^3
=
\sum_{j=1}^{T_n}(2j-1)
=
T_n^2.
\end{equation}

Using the definition of $T_n$, we obtain the modern form of Nicomachus' identity,

\begin{equation}
\boxed{
\sum_{k=1}^{n}k^3
=
\left(\sum_{k=1}^{n}k\right)^2
}
\end{equation}

or equivalently

\begin{equation}
\sum_{k=1}^{n}k^3
=
\left(\frac{n(n+1)}{2}\right)^2.
\end{equation}

The triangular numbers also provide a common modern route to the same result \cite{Bryant1990}. Indeed,

\begin{equation}
T_k^2-T_{k-1}^2=k^3.
\end{equation}

Summing this identity from $k=1$ to $k=n$ gives

\begin{equation}
\begin{aligned}
\sum_{k=1}^{n}k^3
&=\sum_{k=1}^{n}\left(T_k^2-T_{k-1}^2\right)\\
&=(T_1^2-T_0^2)+(T_2^2-T_1^2)+\cdots +(T_n^2-T_{n-1}^2)\\
&=T_n^2,
\end{aligned}
\end{equation}

since all intermediate terms cancel. This proof is short and efficient, but it begins by introducing the triangular numbers and then verifies that their consecutive squares differ by a cube. The historical odd-number pattern suggests the reverse question: why should consecutive triangular numbers appear in the first place? \\The next section approaches that question by starting instead from a difference of squares.

\section{A Difference-of-Squares Route and the Special Role of Cubes}

Nicomachus' observation raises a natural question. Is the fact that $n^3$ can be written as a sum of $n$ consecutive odd integers peculiar to cubes, or does the same phenomenon occur for other powers? To investigate this question, we start from a difference of squares and ask when it represents a block of exactly $n$ consecutive odd integers.

Let $n$ and $p$ be positive integers, with $p\geq 2$, and look for integers $b>a\geq0$ such that

\begin{equation}
n^p=b^2-a^2.
\end{equation}

Before imposing any further condition, it is useful to recall how such difference-of-squares representations are related to factorizations of $n^p$. Set

\begin{equation}
\ell=b-a,
\qquad
m=b+a.
\end{equation}

Then

\begin{equation}
\ell m=n^p.
\end{equation}

Conversely, given a factorization $n^p=\ell m$ with $m\geq\ell>0$ and with $\ell$ and $m$ of the same parity, one recovers

\begin{equation}
b=\frac{m+\ell}{2},
\qquad
a=\frac{m-\ell}{2}.
\end{equation}

Thus different suitable factor pairs of $n^p$ may give different representations of $n^p$ as a difference of two squares. The representation relevant to Nicomachus' pattern will be selected by requiring the difference to contain exactly $n$ consecutive odd integers.

A difference of two squares has a natural interpretation in terms of consecutive odd integers. Indeed,

\begin{equation}
\begin{aligned}
b^2-a^2
&=\bigl((a+1)^2-a^2\bigr)
 +\bigl((a+2)^2-(a+1)^2\bigr)
 +\cdots
 +\bigl(b^2-(b-1)^2\bigr)\\
&=(2a+1)+(2a+3)+\cdots +(2b-1).
\end{aligned}
\end{equation}

The odd integers on the last line may be indexed by

\begin{equation}
a+1,a+2,\ldots,b.
\end{equation}

Hence their number is

\begin{equation}
b-a-1+1=b-a.
\end{equation}

Thus a difference $b^2-a^2$ is the sum of exactly $b-a$ consecutive odd integers. If we want $n^p$ to be represented by exactly $n$ such integers, as in Nicomachus' pattern for $n^3$, we must therefore require

\begin{equation}
b-a=n.
\end{equation}

Factoring the difference of squares gives

\begin{equation}
(b-a)(b+a)=n^p.
\end{equation}

Using $b-a=n$, we obtain

\begin{equation}
b+a=n^{p-1},
\end{equation}

and therefore

\begin{equation}
\boxed{
a=\frac{n^{p-1}-n}{2},
\qquad
b=\frac{n^{p-1}+n}{2}.
}
\end{equation}

These quantities are integers, since $n^{p-1}$ and $n$ have the same parity. Hence, for every $p\geq2$,

\begin{equation}
n^p
=
\left(\frac{n^{p-1}+n}{2}\right)^2
-
\left(\frac{n^{p-1}-n}{2}\right)^2,
\end{equation}

so every $n^p$ can be written as a sum of $n$ consecutive positive odd integers. This particular case, and related representations of powers by consecutive odd integers, have been considered in several forms in the literature \cite{Nair1993,Nelsen1993,Junaidu2010,Vasquez2013}.

Thus the local feature observed by Nicomachus---that $n^3$ is a sum of $n$ consecutive odd integers---is not, by itself, peculiar to cubes. The special feature appears when we ask what the two integers $a$ and $b$ become in the cubic case.

For $p=3$, the formulas above reduce to

\begin{equation}
a=\frac{n^2-n}{2}
=\frac{n(n-1)}{2}
=T_{n-1},
\end{equation}

and

\begin{equation}
b=\frac{n^2+n}{2}
=\frac{n(n+1)}{2}
=T_n.
\end{equation}

Consequently,

\begin{equation}
\boxed{
n^3=T_n^2-T_{n-1}^2.
}
\end{equation}

The triangular numbers that appear in the standard telescoping proof are therefore not introduced in advance: they emerge from the difference-of-squares representation of $n^3$ by $n$ consecutive odd integers.

Moreover, the cubic case is the only power for which this construction produces consecutive triangular numbers for every $n$. Recall that the construction already imposes

\begin{equation}
b-a=n,
\end{equation}

while consecutive triangular numbers satisfy

\begin{equation}
T_n-T_{n-1}=n.
\end{equation}

Therefore it is enough to determine when $b=T_n$. Indeed, if $b=T_n$, then

\begin{equation}
a=b-n=T_n-n=T_{n-1}.
\end{equation}

Using the expression previously obtained for $b$, the condition $b=T_n$ becomes

\begin{equation}
\frac{n^{p-1}+n}{2}=\frac{n^2+n}{2},
\end{equation}

and hence

\begin{equation}
n^{p-1}=n^2.
\end{equation}

For any $n>1$, this forces

\begin{equation}
p=3.
\end{equation}

Thus, among all powers $n^p$ represented as sums of $n$ consecutive odd integers by this construction, cubes are uniquely characterized by the fact that the associated difference of squares is formed from consecutive triangular numbers.

This triangular structure also explains the global pattern recorded by Nicomachus. For $p=3$, the block representing $n^3$ is

\begin{equation}
(2T_{n-1}+1)+(2T_{n-1}+3)+\cdots +(2T_n-1).
\end{equation}

The next block, representing $(n+1)^3$, begins with

\begin{equation}
2T_n+1.
\end{equation}

Thus the successive blocks fit together without gaps or repetitions:

\begin{equation}
\underbrace{1}_{1^3},
\qquad
\underbrace{3+5}_{2^3},
\qquad
\underbrace{7+9+11}_{3^3},
\qquad
\underbrace{13+15+17+19}_{4^3},
\qquad\ldots
\end{equation}

The distinction is therefore not that cubes alone can be expressed as sums of $n$ consecutive odd integers. Rather, only for cubes does this representation recover consecutive triangular numbers, and it is precisely this triangular structure that makes the individual blocks concatenate into Nicomachus' sequence.

It is worth noting that the condition $b-a=n$ was imposed specifically to preserve Nicomachus' $n$-term pattern; it is not the only possible choice. More generally, suppose that we require

\begin{equation}
b-a=n^r,
\end{equation}

so that the corresponding difference of squares contains exactly $n^r$ consecutive odd integers. Since

\begin{equation}
(b-a)(b+a)=n^p,
\end{equation}

we then obtain

\begin{equation}
b+a=n^{p-r},
\end{equation}

and hence

\begin{equation}
\boxed{
a=\frac{n^{p-r}-n^r}{2},
\qquad
b=\frac{n^{p-r}+n^r}{2}.
}
\end{equation}

Therefore, for $1\leq r\leq\lfloor p/2\rfloor$,

\begin{equation}
n^p
=
\left(\frac{n^{p-r}+n^r}{2}\right)^2
-
\left(\frac{n^{p-r}-n^r}{2}\right)^2,
\end{equation}

or equivalently

\begin{equation}
n^p
=
(n^{p-r}-n^r+1)
+(n^{p-r}-n^r+3)
+\cdots
+(n^{p-r}+n^r-1),
\end{equation}

a sum of exactly $n^r$ consecutive odd integers.

Thus the representation of a power as a block of consecutive odd integers need not be unique: different choices of $r$ correspond to different factorizations

\begin{equation}
n^p=n^r n^{p-r},
\end{equation}

and hence to blocks of different lengths. For example,

\begin{equation}
3^5=79+81+83
\end{equation}

uses $3=3^1$ consecutive odd integers, whereas

\begin{equation}
3^5
=
19+21+23+25+27+29+31+33+35
\end{equation}

uses $9=3^2$ consecutive odd integers. Such more general representations are known in the literature \cite{Vasquez2013}.

The choice $r=1$ is the one singled out here because it preserves the defining feature of Nicomachus' pattern: the $n$th power is represented by exactly $n$ consecutive odd integers.

\section{A Faulhaber Perspective: Global Uniqueness}

The previous section showed that the cubic case is exceptional at the level of the individual odd-number blocks: among the representations of $n^p$ as sums of $n$ consecutive odd integers, only $p=3$ produces consecutive triangular numbers. It is natural to ask whether Nicomachus' identity is also exceptional at the level of the corresponding power sums.

For a positive integer $s$, let

\begin{equation}
S_s(N)=\sum_{k=1}^{N}k^s.
\end{equation}

For each fixed positive integer $s$, the power sum $S_s(N)$ is a polynomial
in $N$ of degree $s+1$. These polynomial expressions are commonly referred to as Faulhaber polynomials or Faulhaber formulas \cite{Knuth1993}. For
example,

\begin{equation}
S_1(N)=\frac{N(N+1)}{2},
\end{equation}

while

\begin{equation}
S_2(N)=\frac{N(N+1)(2N+1)}{6}.
\end{equation}

In general, the leading term of $S_s(N)$ is

\begin{equation}
\frac{1}{s+1}N^{s+1}.
\end{equation}

Thus, for the argument below, only the degree and the leading coefficient
of the Faulhaber polynomial are needed.

We ask for which positive integers $s$ and $r$ an identity of the form

\begin{equation}
S_s(N)=S_1(N)^r
\end{equation}

can hold for every positive integer $N$. Faulhaber's formula implies that $S_s(N)$ is a polynomial in $N$ of degree $s+1$ with leading coefficient $1/(s+1)$ \cite{Knuth1993}. On the other hand,

\begin{equation}
S_1(N)^r
=
\left(\frac{N(N+1)}{2}\right)^r
\end{equation}

is a polynomial of degree $2r$ with leading coefficient $1/2^r$.

If the two expressions agree for every positive integer $N$, then the corresponding polynomials are identical. Their degrees and leading coefficients must therefore coincide. Hence

\begin{equation}
s+1=2r
\end{equation}

and

\begin{equation}
\frac{1}{s+1}=\frac{1}{2^r}.
\end{equation}

The second condition gives

\begin{equation}
s+1=2^r.
\end{equation}

Combining the two relations, we obtain

\begin{equation}
2r=2^r.
\end{equation}

For positive integer $r$, the only solutions are $r=1$ and $r=2$. Indeed, both satisfy the equation, while for $r\geq3$ one has $2^r>2r$. Therefore the only possible pairs are

\begin{equation}
(s,r)=(1,1)
\end{equation}

and

\begin{equation}
(s,r)=(3,2).
\end{equation}

The first is the trivial identity $S_1(N)=S_1(N)$. The second is precisely Nicomachus' identity,

\begin{equation}
\boxed{
S_3(N)=S_1(N)^2.
}
\end{equation}

Thus the special role of cubes appears in two complementary ways. Locally, only the cubic case makes the difference-of-squares representation by $n$ consecutive odd integers involve consecutive triangular numbers. Globally, apart from the trivial case, Nicomachus' identity is the unique relation of the form $S_s(N)=S_1(N)^r$ among positive integer exponents.

\section{Conclusion}

Nicomachus' identity is often encountered in the compact form

\begin{equation}
\sum_{k=1}^{n}k^3
=
\left(\sum_{k=1}^{n}k\right)^2,
\end{equation}

but its structure is more transparent when one returns to the original odd-number pattern. Nicomachus' observation organizes successive cubes into successive blocks of consecutive odd integers, and the total number of odd integers appearing in the first $n$ blocks is the triangular number $T_n$. This immediately leads to the familiar identity through the fact that the sum of the first $T_n$ odd integers is $T_n^2$.

The difference-of-squares viewpoint clarifies why triangular numbers arise. Requiring a representation of $n^p$ as a sum of exactly $n$ consecutive odd integers leads to

\begin{equation}
n^p
=
\left(\frac{n^{p-1}+n}{2}\right)^2
-
\left(\frac{n^{p-1}-n}{2}\right)^2.
\end{equation}

This construction exists for every $p\geq2$, so the representation by $n$ consecutive odd integers is not itself peculiar to cubes. What distinguishes the cubic case is that the two integers appearing in the difference of squares become precisely the consecutive triangular numbers $T_{n-1}$ and $T_n$. This is what makes the individual odd-number blocks fit together without gaps or repetitions and recovers Nicomachus' pattern.

The Faulhaber argument gives a complementary global characterization. Among identities of the form

\begin{equation}
S_s(N)=S_1(N)^r,
\end{equation}

with positive integer exponents, the only possibilities are the trivial case $(s,r)=(1,1)$ and Nicomachus' pair $(s,r)=(3,2)$. Thus the cubic identity is exceptional in two related senses: locally, through the triangular structure hidden in its difference-of-squares representation, and globally, through the uniqueness of the corresponding power-sum identity.

\section*{Disclosure Statement}
The author declares that there are no conflicts of interest related to this work.

\end{document}